\documentclass[a4paper,12pt]{article}
\usepackage{comment}
\usepackage{cite}
\usepackage{amsmath}
\usepackage{amssymb}
\usepackage{amsfonts}
\usepackage[T1]{fontenc}
\usepackage[utf8]{inputenc}
\usepackage{graphicx}
\usepackage{fancyhdr}
\usepackage{float}
\usepackage{xcolor}
\usepackage{authblk}
\usepackage{mathrsfs}
\usepackage{empheq}
\usepackage{url}

\usepackage{geometry}
\begin{document}

\title{A simple upper bound for the perimeter of an ellipse}

\author[$\dagger$]{Jean-Christophe {\sc Pain}\\
\small
CEA, DAM, DIF, F-91297 Arpajon, France\\
Université Paris-Saclay, CEA, Laboratoire Matière en Conditions Extr\^emes,\\ 
91680 Bruyères-le-Ch\^atel, France
}

\maketitle

\begin{abstract}
We propose a simple derivation of an upper bound for the perimeter of an ellipse. The procedure, which relies on the use of elliptic integrals, consists in introducing, \emph{via} inequalities and convexity properties, specific integrals which can be calculated analytically. 
\end{abstract}

\section{Introduction}

The perimeter of an ellipse is the length of the continuous line forming its boundary. Computing accurate approximations to the perimeter of an ellipse has been a subject of interest for mathematicians for a long time \cite{Adlaj2012,Borwein1987,Bailey2016}. Many approaches have been used to estimate its value, such as approximation formulas (two formulas derived by Ramanujan are particularly famous \cite{Almkvist1988,Villarino} and we provide a non-exhaustive list in Appendix B) \cite{Hiriart,collect,Barnard2001}, infinite series (Maclaurin, Gauss-Kummer, Euler, etc.), hypergeometric functions \cite{Abbott2009,Chandrupatla2010} (see Appendix C), and of course (elliptic) integrals \cite{Abel1827,Abramowitz1972,Arfken1985,Prudnikov1986,Whittaker1990}.

It is well known that the perimeter of an ellipse can be expressed exactly as a complete elliptic integral of the second kind. Approximate formulas can, of course, be obtained by truncating the series representations of exact formulas. The difficulty of calculating exactly the perimeter has aroused a particular interest in the search for bounds. A long time ago, Kepler found the geometric mean of semi-major and semi-minor axis, as a lower bound for the perimeter, and many sophisticated approaches were developed along the years (see for instance \cite{Pfiefer1988,Wang2012,Gusic2015}). The subject is still an active field of research \cite{He2020,Zhao2022}, important for astrophysical applications (trajectories of planets, satellites, \emph{etc.}).

An ellipse with semi-major axis $a$ and semi-minor axis $b$ is parametrized by $x=a\cos\phi$ and $y=b\sin\phi$. The infinitesimal variation of curvilinear abscissa $s$ is $ds=\sqrt{dx^2+dy^2}$. The length of a quarter of the ellipse is provided by the elliptic integral
\begin{equation}\label{def}
L=\int_0^{\pi/2}ds=\int_0^{\pi/2}\sqrt{a^2\sin^2\phi+b^2\cos^2\phi}~d\phi.
\end{equation}
We will use the latter formula to derive the upper bound. The total perimeter of the ellipse will therefore be
\begin{equation}
L_e=4L.
\end{equation}
The well-known simple lower bound $2\pi\sqrt{ab}$ and upper bound $4a+\pi b$ are recalled in sections \ref{sec1} and \ref{sec2} respectively, and the main result of the present work, which is an upper bound more accurate than the one of section \ref{sec2}, but still rather simple, is presented in section \ref{sec3}.

\section{The geometric-mean lower bound}\label{sec1}

Let us consider two parametrizations of an ellipse
\begin{equation}
\begin{array}{l}
x_1=a\cos\phi,\\
y_1=b\sin\phi
\end{array}
\end{equation}
and
\begin{equation}
\begin{array}{l}
x_2=b\cos\phi,\\
y_2=a\sin\phi.
\end{array}
\end{equation}
The two associated lengths of the quarter of the ellipse are
\begin{equation}
L_1=\int_0^{\pi/2}\sqrt{a^2\cos^2\phi+b^2\sin^2\phi}~d\phi
\end{equation}
and
\begin{equation}
L_2=\int_0^{\pi/2}\sqrt{b^2\cos^2\phi+a^2\sin^2\phi}~d\phi.
\end{equation}
Since $t\rightarrow\sqrt{t}$ is a concave function, one has
\begin{equation}
\sqrt{\lambda t_1+(1-\lambda)t_2}\geq\lambda\sqrt{t_1}+(1-\lambda)\sqrt{t_2}, 
\end{equation}
yielding ($\lambda=\cos^2\phi$, $t_1=a^2$ and $t_2=b^2$):
\begin{equation}
L_1\geq\int_0^{\pi/2}\left(a\cos^2\phi+b\sin^2\phi\right)d\phi
\end{equation}
and
\begin{equation}
L_2\geq\int_0^{\pi/2}\left(b\cos^2\phi+a\sin^2\phi\right)d\phi
\end{equation}
and therefore
\begin{equation}
L=L_1=L_2=\frac{1}{2}\left(L_1+L_2\right)\geq (a+b)\frac{\pi}{4}.
\end{equation}
Setting $ab=R^2$, the minimum of $a+b=a+R/a^2=f(a)$ is obtained for $R=\sqrt{ab}$ and is equal to $2R$. In other words, the ellipse which has the minimum perimeter for a given area is the circle of radius $R$. This yields the well-known result
\begin{align}
\frac{\pi}{2}\sqrt{ab}\leq L
\end{align}
and for the total perimeter of the ellipse
\begin{equation}
2\pi\sqrt{ab}\leq L_e.
\end{equation}

\begin{figure}
\centering
\includegraphics[width=10cm]{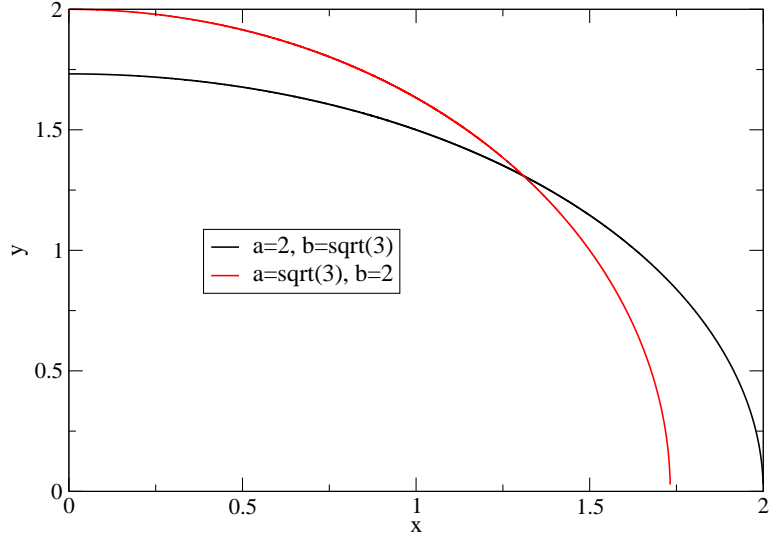}
\caption{quarters of ellipse with semi-major and semi-minor axis $(a=2,b=\sqrt{3})$ and $(a=2,b=\sqrt{3})$ respectively.}
\label{fig23}
\end{figure}

\section{First upper bound}\label{sec2}

Equation (\ref{def}) is equivalent to
\begin{equation}
L=a\int_0^{\pi/2}\sqrt{1-\frac{\left(a^2-b^2\right)}{a^2}\cos^2\phi}~d\phi
\end{equation}
which can be put in the form
\begin{equation}
L=a\frac{\pi}{2}\mathscr{F}\left(\frac{b^2-a^2}{a^2}\right)\cos^2\phi~ d\phi=a\frac{\pi}{2}\mathscr{F}\left(-e^2\right)
\end{equation}
where we have introduced the function
\begin{equation}
\mathscr{F}(t)=\frac{2}{\pi}\int_0^{\pi/2}\sqrt{1+t\cos^2\phi}~d\phi
\end{equation}
where $e$ represents the eccentricity of the ellipse defined by
\begin{equation}
\frac{b^2}{a^2}=1-e^2.
\end{equation}
The asymptotic behaviour of the function $\mathscr{F}$ is studied in Appendix A. If $-1<t<0$, then $-t/(1+t)>0$ and
\begin{equation}
\mathscr{F}(t)-\mathscr{F}(-1)=\frac{2}{\pi}\int_0^{\pi/2}\left(\sqrt{1+t\cos^2\phi}-\sqrt{1-\cos^2\phi}\right)d\phi.
\end{equation}
Since
\begin{equation}
\sqrt{u}-\sqrt{v}=\frac{u-v}{\sqrt{u}+\sqrt{v}},
\end{equation}
we get
\begin{equation}\label{inter}
\mathscr{F}(t)-\mathscr{F}(-1)=\frac{2}{\pi}\int_0^{\pi/2}\frac{(1+t)\cos^2\phi}{\sqrt{1+t\cos^2\phi}+\sqrt{1-\cos^2\phi}}d\phi.
\end{equation}
and
\begin{equation}
\mathscr{F}(t)-\mathscr{F}(-1)\leq\frac{2}{\pi}(1+t)\int_0^{\pi/2}\frac{\cos^2\phi}{\sqrt{1+t}}d\phi,
\end{equation}
\emph{i.e.}
\begin{equation}
\mathscr{F}(t)-\mathscr{F}(-1)\leq\frac{2}{\pi}\sqrt{1+t}~\frac{\pi}{4}.
\end{equation}
because
\begin{equation}
\int_0^{\pi/2}\cos^2\phi~d\phi=\frac{\pi}{4}.
\end{equation}
This enables us to write
\begin{equation}
\mathscr{F}(t)-\mathscr{F}(-1)\leq\frac{\sqrt{1+t}}{2}.
\end{equation}
Since $-1<-e^2<0$, we get
\begin{equation}
\mathscr{F}(-e^2)-\mathscr{F}(-1)\leq\frac{\sqrt{1-e^2}}{2},
\end{equation}
leading to
\begin{equation}\label{firs}
L\leq a+\frac{\pi}{4}b.
\end{equation}
and for the total perimeter of the ellipse
\begin{empheq}[box=\fbox]{align}\label{firstup}
L_e\leq 4a+\pi b.
\end{empheq}

\begin{figure}
    \centering
    \includegraphics[width=10cm]{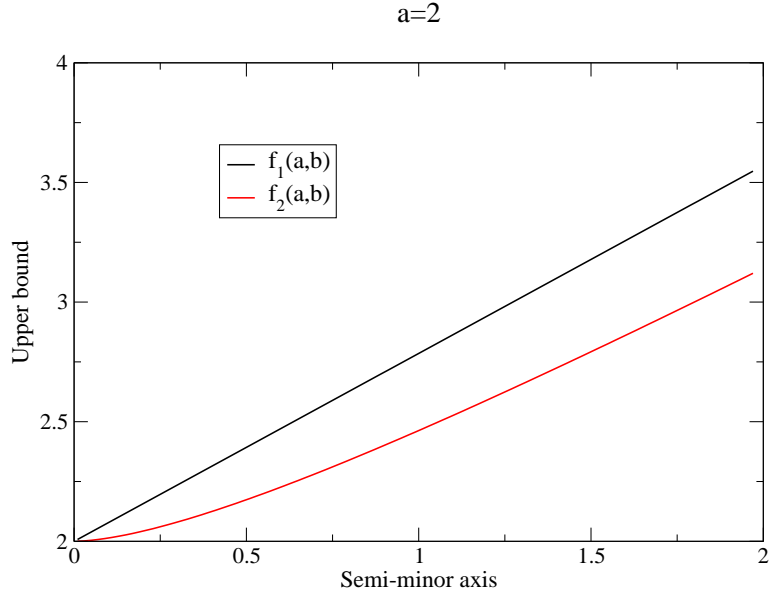}
    \caption{Comparison between the first and second upper bounds of the perimeter of a quarter of an ellipse with semi-major axis $a=2$. On has $f_1(a,b)=a+\frac{\pi}{4}b$ (Eq. (\ref{firs})) and $f_2(a,b)=a+\frac{b^2}{a^2}\left[\sqrt{a^2-b^2}\ln\left(\frac{a+\sqrt{a^2-b^2}}{b}\right)+\frac{\pi}{2}b-a\right]$ (Eq. (\ref{seco})).}
    \label{fig1}
\end{figure}

\section{Second upper bound}\label{sec3}

Let us go back to Eq. (\ref{inter}):
\begin{equation}
\mathscr{F}(t)-\mathscr{F}(-1)=\frac{2}{\pi}\int_0^{\pi/2}\frac{(1+t)\cos^2\phi}{\sqrt{1+t\cos^2\phi}+\sqrt{1-\cos^2\phi}}d\phi.
\end{equation}
We have
\begin{equation}
\mathscr{F}(t)-\mathscr{F}(-1)\leq\frac{2}{\pi}(1+t)\int_0^{\pi/2}\frac{\cos^2\phi}{\sqrt{1+t}+\sin\phi}d\phi,
\end{equation}
where we introduce, for $0<\xi<1$, the following integral
\begin{equation}
\mathscr{K}(\xi)=\int_0^{\pi/2}\frac{\cos^2\phi}{\xi+\sin\phi}d\phi.
\end{equation}
We obtain
\begin{eqnarray}
\mathscr{K}(\xi)&=&\int_0^{\pi/2}\frac{1-\sin^2\phi}{\xi+\sin\phi}d\phi=\int_0^{\pi/2}\frac{\left(\xi^2-\sin^2\phi\right)+\left(1-\xi^2\right)}{1+\sin\phi}d\phi\nonumber\\
&=&\int_0^{\pi/2}\left(\xi-\sin\phi\right)d\phi+\left(1-\xi^2\right)\int_0^{\pi/2}\frac{d\phi}{\xi+\sin\phi}\nonumber\\
&=&\xi\frac{\pi}{2}-1+\left(1-\xi^2\right)\mathscr{J}(\xi)
\end{eqnarray}
with
\begin{equation}
\mathscr{J}(\xi)=\int_0^{\pi/2}\frac{d\phi}{\xi+\sin\phi}.
\end{equation}
Setting $t=\tan(\phi/2)$, we have 
\begin{equation}
dt=\frac{1}{2}\left[1+\tan^2\left(\frac{\phi}{2}\right)\right]d\phi,    
\end{equation}
\emph{i.e.} $d\phi=2dt/(1+t^2)$. Then, since 
\begin{equation}
\sin\phi=\frac{2t}{\left(1+t^2\right)},    
\end{equation}
we get
\begin{equation}
\frac{d\phi}{\xi+\sin\phi}=\frac{2dt}{\xi t^2+2t+\xi}.
\end{equation}
Writing
\begin{equation}
\frac{2}{\xi t^2+2t+\xi}=\frac{2}{\xi(t-\chi_1)(t-\chi_2)}
\end{equation}
with
\begin{equation}
\chi_1=\frac{-1-\sqrt{1-\xi^2}}{\lambda}
\end{equation}
and
\begin{equation}
\chi_2=\frac{-1+\sqrt{1-\xi^2}}{\lambda},
\end{equation}
enables us to write
\begin{eqnarray}
\mathscr{J}(\xi)&=&\frac{2}{\xi}\int_0^1\frac{dt}{(t-\chi_1)(t-\chi_2)}=\frac{2}{\xi(\chi_1-\chi_2)}\left[\frac{1}{t-\chi_1}-\frac{1}{t-\chi_2}\right]dt\nonumber\\
&=&\frac{2}{\xi(\chi_1-\chi_2)}\left|\ln\left(\frac{t-\chi_1}{t-\chi_2}\right)\right|_0^1=\frac{2}{\xi(\chi_1-\chi_2)}\ln\left[\frac{(1-\chi_1)\chi_2}{(1-\chi_2)\chi_1}\right].
\end{eqnarray}
We have also
\begin{equation}
\chi_1\chi_2=1
\end{equation}
and
\begin{equation}
\chi_1-\chi_2=-\frac{2}{\xi}\sqrt{1-\xi^2}
\end{equation}
and thus
\begin{eqnarray}
\mathscr{J}(\xi)&=&-\frac{1}{\sqrt{1-\xi^2}}\ln\left(\frac{\chi_2-1}{\chi_1-1}\right)=\frac{1}{\sqrt{1-\xi^2}}\ln\left(\frac{\xi+1+\sqrt{1-\xi^2}}{\xi+1-\sqrt{1-\xi^2}}\right)\nonumber\\
&=&\frac{1}{\sqrt{1-\xi^2}}\ln\left(\frac{\sqrt{\xi+1}+\sqrt{1-\xi}}{\sqrt{\xi+1}-\sqrt{1-\xi}}\right)\nonumber\\
&=&\frac{1}{\sqrt{1-\xi^2}}\ln\left(\frac{1+\sqrt{1-\xi^2}}{\xi}\right),
\end{eqnarray}
yielding
\begin{equation}
\mathscr{K}(\xi)=\frac{\xi\pi}{2}-1+\sqrt{1-\xi^2}~\ln\left(\frac{1+\sqrt{1-\xi^2}}{\xi}\right).
\end{equation}
Therefore
\begin{equation}
L=\frac{\pi}{2}a\mathscr{F}\left(\frac{b^2}{a^2}-1\right)\leq\frac{\pi}{2}a\left[\mathscr{F}(-1)+\frac{2}{\pi}\left(\frac{b^2}{a^2}\right)\mathscr{K}\left(\frac{b}{a}\right)\right],
\end{equation}
or finally
\begin{equation}\label{seco}
L\leq a+\frac{b^2}{a^2}\left[\sqrt{a^2-b^2}\ln\left(\frac{a+\sqrt{a^2-b^2}}{b}\right)+\frac{\pi}{2}b-a\right]
\end{equation}
and for the total perimeter of the ellipse
\begin{empheq}[box=\fbox]{align}
L_e\leq 4a+4\frac{b^2}{a^2}\left[\sqrt{a^2-b^2}~\ln\left(\frac{a+\sqrt{a^2-b^2}}{b}\right)+\frac{\pi}{2}b-a\right].
\end{empheq}
This upper bound is better than the first one in Eq. (\ref{firstup}) (see figure \ref{fig1}).

\section{Conclusion}

In this article, we obtained an upper bound for the perimeter of an ellipse. The calculation consists in simplifying elliptic integrals by applying inequalities and convexity properties. It is hoped that the method may stimulate the derivation of further bounds, improving the compromise between simplicity and accuracy.

\section*{Appendix A: Asymptotic behavior}

Let us consider $t>-1$, i.e. $-t/(1+t)>-1$. We have
\begin{eqnarray}
\mathscr{F}\left(\frac{-t}{1+t}\right)&=&\frac{2}{\pi}\int_0^{\pi/2}\frac{\sqrt{1+t-t\cos^2\phi}}{\sqrt{1+t}}d\phi\nonumber\\
&=&\frac{2}{\pi}\frac{1}{\sqrt{1+t}}\int_0^{\pi/2}\sqrt{1+t\sin^2\phi}~d\phi
\end{eqnarray}
Making the change of variables $\psi=\pi/2-\phi$, we get
\begin{equation}
\int_0^{\pi/2}\sqrt{1+t\sin^2\phi}~d\phi=\int_0^{\pi/2}\sqrt{1+t\cos^2\phi}~d\phi
\end{equation}
and thus
\begin{equation}
\sqrt{1+t}~\mathscr{F}\left(-\frac{t}{1+t}\right)=\frac{2}{\pi}\int_0^{\pi/2}\sqrt{1+t\cos^2\phi}~d\phi
\end{equation}
and
\begin{equation}
\mathscr{F}(t)=\sqrt{1+t}~\mathscr{F}\left(-\frac{t}{1+t}\right)
\end{equation}
\begin{equation}
\frac{\mathscr{F}(t)}{\sqrt{t}}=\frac{\mathscr{F}(t)}{\sqrt{1+t}}\frac{\sqrt{1+t}}{\sqrt{t}}=\mathscr{F}\left(-\frac{t}{1+t}\right)\sqrt{1+\frac{1}{t}}
\end{equation}
If $t\rightarrow\infty$, $-t/(1+t)\rightarrow -1$ and $\mathscr{F}(t)/\sqrt{t}\rightarrow\mathscr{F}(-1)=2/\pi$. Therefore
\begin{equation}
\mathscr{F}(t)\approx\frac{2}{\pi}\sqrt{t}.
\end{equation}

\section*{Appendix B: A few approximate formulas}

In this appendix, we provide a non-exhaustive of approximate formulas for the perimeter of an ellipse:

\vspace{1cm}
\noindent Kepler (1609) [ lower bound of section \ref{sec1} ]:
\begin{equation}
L_e\approx 2\pi\sqrt{ab}
\end{equation}
\begin{equation}
L_e\approx \pi(a+b)
\end{equation}
Euler
\begin{equation}
L_e\approx 2\pi\sqrt{\frac{a^2+b^2}{2}}
\end{equation}
Sipos
\begin{equation}
L_e\approx 2\pi\sqrt{\frac{(a+b)^2}{(\sqrt{a}+\sqrt{b})^2}}
\end{equation}
Cesaro
\begin{equation}
L_e\approx \pi(a+b)+\frac{\pi}{4}\frac{(a-b)^2}{a+b}
\end{equation}
Muir
\begin{equation}
L_e\approx 2\pi\left[\frac{a^{3/2}+b^{3/2}}{2}\right]^{2/3}
\end{equation}
Peano, Boussinesq
\begin{equation}
L_e\approx \pi\left[\frac{3}{2}(a+b)-\sqrt{ab}\right]
\end{equation}
Almkvist
\begin{equation}
L_e\approx 2\pi\frac{2(a+b)^2-\left(\sqrt{a}-\sqrt{b}\right)^4}{\left(\sqrt{a}+\sqrt{b}\right)^2+2\sqrt{2(a+b)}~(ab)^{1/4}}
\end{equation}
Lindner
\begin{equation}
L_e\approx \pi(a+b)\left[1+\frac{1}{8}\left(\frac{a-b}{a+b}\right)^2\right]^2
\end{equation}
Ramanujan: first formula
\begin{equation}
L_e\approx\pi\left[3(a+b)-\sqrt{(3a+b)(a+3b)}\right]    
\end{equation}
Ramanujan: second formula
\begin{equation}
L_e\approx\pi\left[(a+b)-\frac{3(a-b)^2}{10(a+b)+\sqrt{a^2+14ab+b^2}}\right].
\end{equation}
\section*{Appendix C: Expressions in terms of special functions or infinite sums}

The integral $\mathscr{F}(t)$ can be expanded as
\begin{eqnarray}
\mathscr{F}(t)&=&\frac{2}{\pi}\int_0^{\pi/2}\sqrt{1+t\cos^2\phi}~d\phi\nonumber\\
&=&1+\sum_{k=1}^{\infty}(-1)^{k+1}t^k\frac{(2k-3)!!}{2^kk!}\frac{2}{\pi}\int_0^{\pi/2}\cos^{2k}\phi~d\phi
\end{eqnarray}
with $(2k-3)!!=(2k-3)\cdots 5\times 3\times 1$ and $(-1)!!=1$. Setting
\begin{eqnarray}
\mathscr{A}_k=\int_0^{\pi/2}\cos^{2k}\phi d\phi&=&\int_0^{\pi/2}\cos^{2k-1}\phi\cos\phi~d\phi\nonumber\\
&=&\left|\cos^{2k-1}\phi\sin\phi\right|_0^{\pi/2}+\int_0^{\pi/2}(2k-1)\cos^{2k-2}\phi\sin^2\phi~d\phi\nonumber\\
&=&(2k-1)\mathscr{A}_{k-1}-(2k-1)\mathscr{A}_k,
\end{eqnarray}
we have
\begin{equation}
\mathscr{A}_k=\frac{(2k-1)}{2k}\mathscr{A}_{k-1}=\frac{(2k-1)!!}{2^kk!}\mathscr{A}_0,
\end{equation}
where $\mathscr{A}_0=\pi/2$ and finally
\begin{eqnarray}
\mathscr{F}(t)&=&1+\sum_{k=1}^{\infty}(-1)^{k+1}\left[\frac{(2k-1)!!}{2^kk!}\right]^2\frac{t^k}{2k-1}
\end{eqnarray}
giving the final Euler-Mac-Laurin series
\begin{eqnarray}
L_e=2a\pi\left\{1+\sum_{k=1}^{\infty}(-1)^{k+1}\left[\frac{(2k-1)!!}{2^kk!}\right]^2\frac{\left(\frac{b^2}{a^2}-1\right)^k}{2k-1}\right\}.
\end{eqnarray}
We mention below a few other formulas:
\vspace{1cm}

\noindent McLaurin
\begin{equation}
L_e=2a\pi~_2F_1\left[\begin{array}{c}
-\frac{1}{2}, \frac{1}{2}\\
1
\end{array};1-\frac{b^2}{a^2}\right]
\end{equation}
Gauss-Kummer
\begin{equation}
L_e=\pi(a+b)~_2F_1\left[\begin{array}{c}
-\frac{1}{2}, \frac{1}{2}\\
1
\end{array};\left(\frac{a-b}{a+b}\right)^2\right]
\end{equation}
Cayley
\begin{eqnarray}
L_e&=&4a\left[1+\frac{1}{2}\left(\ln\left(\frac{4a}{b}\right)-\frac{1}{1\times 2}\right)\left(\frac{b}{a}\right)^2+\frac{1^2\times 3}{2^2\times 4}\left(\ln\left(\frac{4a}{b}\right)-\frac{2}{1\times 2}-\frac{1}{3\times 4}\right)\left(\frac{b}{a}\right)^4\right.\nonumber\\
& &\left.+\frac{1^2\times 3^2\times 5}{2^2\times 4^2\times 6}\left(\ln\left(\frac{4a}{b}\right)-\frac{2}{1.2}-\frac{2}{3.4}-\frac{1}{5.6}\right)\left(\frac{b}{a}\right)^6+\cdots\right]
\end{eqnarray}
Euler
\begin{equation}
L_e=\pi\sqrt{2(a^2+b^2)}~_2F_1\left[\begin{array}{c}
\frac{1}{4}, -\frac{1}{4}\\
1
\end{array};\left(\frac{a^2-b^2}{a^2+b^2}\right)^2\right]
\end{equation}
Abbott
\begin{equation}
L_e=2\pi\sqrt{ab}~\mathscr{P}_{1/2}\left(\frac{a^2+b^2}{2ab}\right),
\end{equation}
where $\mathscr{P}_{1/2}$ is a Legendre function.

\end{document}